\documentclass[10pt]{amsart}
\usepackage{geometry}
\usepackage{amsmath}
\usepackage{amsthm}
\usepackage{amssymb}
\usepackage{mathdots} 
\usepackage{graphicx}
\usepackage[multiple]{footmisc}
\usepackage{multicol}
\usepackage{relsize}
\usepackage{tikz}
\usetikzlibrary{arrows}

\usepackage[pdfauthor={Kameryn J Williams},
    pdftitle={The omega-th inner mantle},
    hidelinks 
]{hyperref}

\theoremstyle{plain}
\ifdefined\theorem \else \newtheorem{theorem}{Theorem} \fi
\ifdefined\maintheorem \else \newtheorem{maintheorem}[theorem]{Main Theorem} \fi
\ifdefined\theoremschema \else \newtheorem{theoremschema}[theorem]{Theorem Schema} \fi
\ifdefined\lemma \else \newtheorem{lemma}[theorem]{Lemma} \fi
\ifdefined\lemmaschema \else \newtheorem{lemmaschema}[theorem]{Lemma Schema} \fi
\ifdefined\proposition \else \newtheorem{proposition}[theorem]{Proposition} \fi
\ifdefined\corollary \else \newtheorem{corollary}[theorem]{Corollary} \fi
\ifdefined\fact \else \newtheorem{fact}[theorem]{Fact} \fi
\ifdefined\problem \else \newtheorem{problem}[theorem]{Problem} \fi
\ifdefined\conjecture \else \newtheorem{conjecture}[theorem]{Conjecture} \fi
\ifdefined\question \else \newtheorem{question}[theorem]{Question} \fi
\ifdefined\observation \else \newtheorem{observation}[theorem]{Observation} \fi
\newtheorem*{theorem*}{Theorem}
\newtheorem*{lemma*}{Lemma}
\newtheorem*{theoremschema*}{Theorem Schema}
\newtheorem*{proposition*}{Proposition}
\newtheorem*{conjecture*}{Conjecture}
\newtheorem*{question*}{Question}
\newtheorem*{definition*}{Definition}
\newtheorem{sublemma}{Lemma}[theorem]
\theoremstyle{definition}
\ifdefined\definition \else \newtheorem{definition}[theorem]{Definition} \fi
\ifdefined\subdefinition \else \newtheorem{subdefinition}{Definition}[theorem] \fi
\ifdefined\maindefinition \else \newtheorem{maindefinition}[theorem]{Main Definition} \fi
\theoremstyle{remark}
\ifdefined\remark \else \newtheorem{remark}[theorem]{Remark} \fi
\ifdefined\remarks \else \newtheorem{remarks}[theorem]{Remarks} \fi
\ifdefined\example \else \newtheorem{example}[theorem]{Example} \fi
\ifdefined\claim \else \newtheorem{claim}[theorem]{Claim} \fi
\ifdefined\exercise \else \newtheorem{exercise}[theorem]{Exercise} \fi
\newtheorem*{remark*}{Remark}
\ifdefined\acknowledgment \else \newtheorem*{acknowledgment}{Acknowledgment} \fi
\ifdefined\dedication \else \newtheorem*{dedication}{Dedicated} \fi
\ifdefined\case \else \newtheorem*{case}{Case} \fi

\newcounter{my_enumerate_counter}

\newcommand\comment[1]{}

\newcommand\Lrm{\mathrm{L}}
\newcommand\Mrm{\mathrm{M}}

\newcommand\Vrm{\mathrm{V}}

\newcommand\Pcal{\mathcal{P}}

\newcommand\Mbb{\mathbb{M}}

\newcommand\Obb{\mathbb{O}}
\newcommand\Pbb{\mathbb{P}}
\newcommand\Qbb{\mathbb{Q}}

\newcommand\Sbb{\mathbb{S}}
\newcommand\Tbb{\mathbb{T}}

\newcommand\Zbb{\mathbb{Z}}

\newcommand{\dom}{\operatorname{dom}}

\ifdefined\alt \else
  \newcommand{\alt}{\operatorname{alt}}
\fi

\newcommand{\forces}{\Vdash}

\newcommand\cat{{}^\smallfrown}

\newcommand\axiom{\mathsf}

\newcommand\GCH{\axiom{GCH}}

\newcommand\ZF{\axiom{ZF}}
\newcommand\ZFC{\axiom{ZFC}}

\newcommand\CC{\axiom{CC}}

\newcommand\AC{\axiom{AC}}

\newcommand\class{\mathrm}

\newcommand\HOD{\class{HOD}}
\newcommand\Ord{\class{Ord}}
\newcommand\Card{\class{Card}}

\newcommand\Succ{\class{Succ}}

\newcommand\Add{\class{Add}}

\newcommand{\seq}[1]{{\left\langle #1 \right\rangle}}

\renewcommand{\epsilon}{\varepsilon}

\newcommand\card[1]{\left\lvert #1 \right\rvert}
\newcommand\length{\operatorname{len}}
\newcommand\rest{\upharpoonright}

\newcommand\powerset{\Pcal}

\ifdefined\Form \else
  \newcommand\Form{\axiom{Form}}
\fi

\renewcommand{\phi}{\varphi}

\renewcommand\CC{\mathbb{CC}}
\newcommand\RCC{\mathbb{RCC}}

\newcommand\headsf{\mathsf{head}}
\newcommand\tailsf{\mathsf{tail}}

\title{The $\omega$-th inner mantle}

\DeclareRobustCommand{\okina}{%
  \raisebox{\dimexpr\fontcharht\font`A-\height}{%
    \scalebox{0.8}{`}%
  }%
}
\author{Kameryn J. Williams}
\address[Kameryn J. Williams]{
University of Hawai\okina{}i at M\=anoa \\
Department of Mathematics \\
2565 McCarthy Mall, Keller 401A \\
Honolulu, HI  96822\\
USA}
\email{kamerynw@hawaii.edu}
\urladdr{http://kamerynjw.net} 

\date{\today}

\begin{document}

\begin{abstract}
This article investigates pathological behavior at the first limit stage in the sequence of inner mantles, obtained by iterating the definition of the mantle to get smaller and smaller inner models. I show: $(A)$ it is possible that the $\omega$-th inner mantle is not a definable class; and $(B)$ it is possible that the $\omega$-th inner mantle is a definable class but does not satisfy $\AC$. This answers a pair of questions of Fuchs, Hamkins, and Reitz \cite{fuchs-hamkins-reitz2015}.
\end{abstract}

\maketitle

\section{Introduction}

An inner model $W \subseteq \Vrm$ is a \emph{ground} if $\Vrm$ is a set forcing extension of $W$. More precisely, $W$ is a ground if there is a poset $\Pbb \in W$ and $G \in \Vrm$ which is generic over $W$ for $\Pbb$ so that $\Vrm = W[G]$. Laver \cite{laver2007}
and independently Woodin \cite{woodin2011b}
proved that the grounds are uniformly definable in $\ZFC$, allowing for them to be quantified over in a first-order manner.
The \emph{mantle} $\Mrm$, first defined and studied in \cite{fuchs-hamkins-reitz2015} is the intersection of the grounds. It follows from their work combined with Usuba's proof of the strong downward directed grounds hypothesis \cite{usuba2017} that the mantle is an inner model of $\ZFC$ and that it is the largest set forcing-invariant inner model. Fuchs, Hamkins, and Reitz \cite{fuchs-hamkins-reitz2015} produced a class forcing notion which forces the ground model to be the mantle of the extension. In particular, this implies that the mantle is not absolute; it is consistent that $\Mrm^\Mrm \ne \Mrm$. As such, it is sensible to iterate taking the mantle to get a sequence of smaller and smaller inner models. 

\begin{definition}
The sequence of \emph{inner mantles} is inductively defined as follows. The zeroth inner mantle $\Mrm^0$ is $\Vrm$. Given the $\eta$-th inner mantle $\Mrm^\eta$ the $(\eta+1)$-th inner mantle is $\Mrm^{\eta + 1} = \Mrm^{\Mrm^\eta}$. And for limit stages $\gamma$ we define the $\gamma$-th inner mantle to be $\Mrm^\gamma = \bigcap_{\eta < \gamma} \Mrm^\eta$. If $\eta$ is least so that $\Mrm^\eta = \Mrm^{\eta + 1}$ say that the sequence of inner mantles has \emph{length} $\eta$ or, synonymously, that the sequence of inner mantles \emph{stabilizes} at $\eta$. 
\end{definition}

In previous joint work with Reitz \cite{reitz-williams2019} we showed that there are class forcing notions $\Mbb(\eta)$, uniformly definable in $\eta$, so that forcing with $\Mbb(\eta)$ makes the ground model the $\eta$-th inner mantle of the extension. This answered a question from \cite{fuchs-hamkins-reitz2015}. In the present article I address two remaining questions from \cite{fuchs-hamkins-reitz2015}. 

The inductive definition of the sequence of inner mantles takes place in the meta-theory. Repeatedly taking the mantle of the mantle of $\ldots$ gives increasingly complex definitions, and it is far from clear that we can uniformly define $\Mrm^n$ for $n < \omega$ to thereby obtain that $\Mrm^\omega$ is definable. Fuchs, Hamkins, and Reitz asked whether there were a model of $\ZFC$ whose $\omega$-th mantle is not definable. The first main theorem of this article answers this question in the positive.

\begin{maintheorem} \label{main:definable}
There is a model $N$ of $\ZFC$ so that $(\Mrm^\omega)^N$ is not a definable class in $N$, if there is any model of $\ZFC$. If there are transitive models of $\ZFC$ then we may take $N$ to be transitive. 
\end{maintheorem}

Observe that if $\Mrm^\omega$ is not a definable class then the sequence of inner mantles cannot stabilize before $\omega$.

As remarked above, the mantle is an inner model of $\ZFC$ and thus if $\Mrm^\eta \models \ZFC$ then so does $M^{\eta+1}$. It is straightforward to see that limit stage inner mantles, provided they are definable classes, must satisfy $\ZF$.

\begin{lemma} \label{lem:zf-is-free}
Suppose $\seq{W_i : i < \lambda \in \Ord}$ is a definable decreasing sequence of inner models of $\ZF$, where the tail sequence $\seq{W_j : i < j < \lambda}$ is definable in each $W_i$. Then $W = \bigcap_{i  < \lambda} W_i$ is an inner model of $\ZF$. In particular, if $\gamma$ is a limit ordinal and $\Mrm^\gamma$ is a definable class so that $\Mrm^\eta \models \ZFC$ for all $\eta < \gamma$ then $\Mrm^\gamma$ is an inner model of $\ZF$.
\end{lemma}

\begin{proof}
It is immediate that $W$ is transitive, contains all ordinals, and is closed under the G\"odel operations. It is easy to check that $W$ is almost universal, meaning that if $A \subseteq W$ is a set then there is $B \in W$ which covers $A$. Namely, let $B = \Vrm_\xi \cap W$ where $A \subseteq \Vrm_\xi$. This $B$ is definable because the sequence $\seq{W_i : i \in I}$ is definable. And $B \in W_i$ for all $i < \lambda$ because the tail of the sequence is definable in $W_i$.
This establishes that $W$ is an inner model of $\ZF$, as desired
\end{proof}

Is this result the best possible? Can we prove that limit stage inner mantles must satisfy Choice? The second main theorem of this article answers this question negatively.

\begin{maintheorem} \label{main:choice}
There is a model $N$ of $\ZFC$ whose $\omega$-th inner mantle is definable, but $(\Mrm^\omega)^N \not \models \AC$, assuming there is any model of $\ZFC$. If transitive models of $\ZFC$ exist, then $N$ may be taken to be transitive. 
\end{maintheorem}

The reader should compare these results to analogous results about the iterated $\HOD$ sequence. Similar to the sequence of inner mantles, this sequence is obtained by iterating the definition of $\HOD$.  Harrington (in unpublished work) established the $\HOD^\omega$ analogue of Main Theorem~\ref{main:definable} and McAloon \cite{mcaloon1974} established the $\HOD^\omega$ analogue of Main Theorems~\ref{main:choice}. The reader who is interested in these classical results is encouraged to consult the treatment by Zadro\.zny \cite{zadrozny1983}.
Indeed, Zadro\.zny's article was a key point of inspiration for the current article, and the techniques herein are an adaptation of the methods from the $\HOD$ context to the mantle context.

This article is structured as follow. First is a section surveying some methods for coding sets in inner mantles, which will be used to establish the main theorems. Next is a section establishing Main Theorem~\ref{main:choice}, followed by a section establishing Main Theorem~\ref{main:definable}. These last two sections are independent, and the reader who is only interested in one of the two main results may freely skip the other. I end with a very short section conjecturing how the work in this paper might be improved.

\section{Coding sets into inner mantles and tree iterations} \label{sec:coding}

I this section I discuss the basic method for ensuring that a set appears in the desired inner mantles. This will be the foundation for the main results of this article. 

The models considered in this article will be extensions of $\Lrm$ by adding Cohen sets.
The most basic coding one might use is what one might call \emph{simple Cohen coding}, where a set of ordinals $x$ is coded by the pattern of which cardinals have a subset which is Cohen generic over $\Lrm$. Namely, consider a set of ordinals $x$ without a maximum.\footnote{If $\sup x \in x$, then $x$ is definable from $x \cup \{\sup x + \omega + n : n \in \omega \}$, so this assumption is harmless for our purposes.} 
Say that $x$ is \emph{simple Cohen coded at $\alpha$} if for all $i < \lambda = \sup x$ we have that $\alpha^{+i}$ contains a subset Cohen generic over $\Lrm$ if and only if $i \in x$. Suppose we start in a ground model $\Vrm$ where the interval $[\alpha, \alpha^{+\lambda})$ is \emph{clean for coding}, meaning that no cardinal in the interval contains a subset Cohen generic over $\Lrm$. Then we can force $x$ to be Cohen coded at $\alpha$ via the full support product of $\Add^\Lrm(\alpha^{+i},1)$ for $i \in x$. 

To ensure that we do not accidentally code too much\footnote{Here is a simple example to illustrate the problem if we do not do this. Suppose we first add a Cohen real $c$ then force over $\Vrm[c]$ to add a Cohen generic subset $d$ of $\kappa > \omega$. If $d \in W \subseteq \Vrm$, then $c \in W$, because by density $c$ is coded into $d$.}
when we do $\Add$ forcings in this paper they will always be as defined in $\Lrm$, even if we are forcing over an extension of $\Lrm$. I wish to emphasize, however, that while each iterand is a forcing from $\Lrm$, the entire iteration will not be from $\Lrm$. This suffices for our coding purposes, because we code things by the \emph{locations} of Cohen generics, not by the generics themselves. In this article, all models will be extensions of $\Lrm$ in which all cardinals are preserved.  In this context, $\Add^L(\alpha,1)$ is always $\mathord{<}\alpha$-distributive, as can be seen by an easy covering argument lifting the distributivity from $\Lrm$ to $\Vrm$.

For later arguments it will be convenient to be able to make the move from knowing that the individual Cohen generics used for coding show up in an inner model to knowing that the entire sequence of Cohens is definable in the inner model. For this purpose we make use of a more complicated coding method, which I will call \emph{Cohen coding}. This coding method will use an iteration of simple Cohen coding to ensure the generics are self-encoded. Let me begin by describing the Cohen coding forcing $\CC(x)$ which codes a set $x$ of ordinals, with $\lambda = \sup x \not \in x$, starting at a cardinal $\alpha$. This coding is done via a full-support length $\omega$ iteration. 
\begin{itemize}
\item At stage $0$, do simple Cohen coding starting at $\alpha^+$ to code $x$. This produces a $\lambda$-sequence $\seq{c_i : i \in x}$ of Cohens where $c_i \subseteq \alpha^{+i+1}$. Using an absolute pairing function we may think of the sequence $\seq{c_i : i \in x}$ as a single set $x_1 \subseteq \alpha^{+\alpha}$ of ordinals. Set $\alpha_1 = \alpha^{+\alpha}$.
\item At stage $n+1$, we want to code the set $x_n \subseteq \alpha_n$. Do simple Cohen coding at ${\alpha_n}^+$ to code $x_n$, producing a sequence $\seq{c_i : i \in x_n}$ of Cohens where $c_i \subseteq {\alpha_n}^{+i+1}$. Again using an absolute pairing function we may think of this sequence as a single set $x_{n+1} \subseteq \alpha_{n+1}$, where $\alpha_{n+1} = {\alpha_n}^{+\alpha_n}$.
\item After $\omega$ many stages, we have the full generic for for Cohen coding, which we may think of either as a sequence $\seq{x_n : n < \omega}$ of sets of ordinals or as a sequence $\bar c$ of Cohens. The domain of $\bar c$ is unbounded in the interval $[\alpha,\alpha_\omega)$, where $\alpha_\omega = \sup_n \alpha_n$.
\end{itemize}

Of course, we do not know which cardinals in stage $n+1$ we will be adding Cohens to until we have the generic for stage $n$. So in the ground model we cannot yet identify this. Nonetheless, we can compute the cardinals $\alpha_n$ in the ground model, and thereby determine the interval $[\alpha,\alpha_\omega)$ used for this coding. If not explicitly addressed, I will implicitly assume that the interval used is clean for coding when talking about Cohen coding.

The following lemma gives the self-encoding property of Cohen coding.

\begin{lemma}[Self-encoding lemma] \label{lem:se}
Let $\Vrm[\bar c]$ be a forcing extension in which a set $x$ of ordinals has been Cohen coded in the interval $[\alpha,\alpha_\omega)$. If $W \subseteq \Vrm[\bar c]$ is an inner model which contains each individual Cohen $c_i$ from $\bar c$, then $W$ contains the entire sequence. More generally, suppose $\bar\Vrm$ is a model in which $\Ord$ many sets of ordinals have been Cohen coded in different disjoint intervals, where the starting points $\alpha$ form a definable class. (Possibly other coding was done to produce $\bar\Vrm$, but we assume no other Cohens were added to the coding regions.) If $W \subseteq \bar\Vrm$ is an inner model which contains each individual Cohen from each of these codings, then the entire $\Ord$-length sequence of all Cohens across all codings is definable over $W$.
\end{lemma}

\begin{proof}
It suffices to prove the more general statement. Consider such a $W$ and we want to show that $\bar C$, the $\Ord$-length sequence of all the Cohens is definable over $W$. First, observe that we can identify $\dom \bar C$ simply by looking at which cardinals in the coding regions contain Cohen subsets. So we only have to uniformly identify which Cohen subset of a given cardinal is the one in $\bar C$. Fix a starting point $\alpha$ and look within the coding region $[\alpha^+,\alpha_\omega)$. To identify the Cohens in $\bar C$ in this region it suffices to identify each $x_n$ for this region. We can identify $x_n$ by looking in the region $[\alpha_n,\alpha_{n+1})$. Namely, $i \in x_n$ if and only if there is a Cohen subset of ${\alpha_n}^{+i+1}$. This definition is uniform across all $n$ and all starting points $\alpha$, so we have defined $\bar C$ in $W$.
\end{proof}

If we would like to code $x$ into the mantle, then it is not enough to Cohen code it once, as the following simple example illustrates. Start with $\Lrm$ and force to add a Cohen real $x \subseteq \omega$. Follow this by forcing over $\Lrm[x]$ to Cohen code $x$ at $\omega_1$ to get $\Lrm[x][\bar c]$. Then $\Lrm$ is a ground of $\Lrm[x][\bar c]$ and so $\Mrm^{\Lrm[x][\bar c]} = \Lrm$ does not contain $x$. The fix to this is to Cohen code $x$ cofinally often.

\begin{lemma} \label{lem:tail-agreement}
Let $W$ be a ground of $\Vrm$. Then $W$ and $\Vrm$ agree on a tail about which cardinals contain subsets which are Cohen generic over $\Lrm$.
\end{lemma}

\begin{proof}
Suppose that $\bar C$ is a proper class-sized sequence of Cohen sets which are in $\Vrm$ but not $W$. Then $W \subseteq W[\bar C] \subseteq \Vrm$. So by the intermediate model theorem we get that $W[\bar C]$ is an extension of $W$ by \emph{set} forcing. But set forcing cannot add a proper class of Cohen sets. Contradiction achieved.
\end{proof}

As a consequence of this lemma, if $x$ is Cohen coded cofinally often then $x$ is in the mantle. This can be achieved by an $\Ord$-length product, repeatedly Cohen coding $x$ in different coding regions. For coding sets into deeper inner mantles, it will be convenient to not use contiguous blocks of cardinals when Cohen coding sets. Let me first make an auxiliary definition.

\begin{definition}
If $R$ is a class of cardinals and $\alpha$ and $i$ are ordinals then let $\Succ^i_R(\alpha)$ be the $i$-th successor of $\alpha$ in $R$. More precisely, $\beta = \Succ^i_R(\alpha)$ is the unique cardinal in $R$ so that $R \cap [\alpha,\beta)$ has ordertype $i$. In particular, $\Succ^i_R(0)$ is always the $i$-th cardinal in $R$.\footnote{Note that if $R = \Card$ then $\Succ^i_R(\alpha) = \alpha^{+i}$.}
\end{definition}

It is straightforward to modify the above-described coding apparatus to go along an arbitrary class $R$ rather than along the class of all cardinals. Namely, replace all references to $\alpha^{+i}$ with $\Succ^i_R(\alpha)$. I will use \emph{simple Cohen coding} and \emph{Cohen coding} to refer to the general case of coding along an arbitrary coding region, and often omit explicit reference to the coding region used. We need, of course, that $R$ remains a definable class in the extension if we want to be able to decode after forcing. In this article, the coding regions we consider will be $\Delta_0$- or $\Delta_1$-definable, and so they will be definable in the extension by the same definition. 

Suppose we want to code $x$ into a finite stage inner mantle $\Mrm^k$. Then we need to code $x$ cofinally often, and we need to code those codes cofinally often, and so on. As such, the pattern our coding will follow is that of a tree, namely ${}^{\le k}\Ord$. And we will need more complicated tree patterns for later forcings, so let us first take a step back and consider a general framework.

In this article we will work only with well-founded trees. My trees grow upward, so $s < t$ in $T$ means that $t$ is an extension of $s$. I will have occasion to talk about subtrees whose root is a non-root node in $T$, so I will use the term \emph{rooted subtree} to refer to those subtrees which have the same root as $T$. All supports will be set-support. I first give a formal definition, following Groszek and Jech's more general treatment of iterations along a general partial order \cite{groszek-jech1991},\footnote{There is a minor implementation difference; Groszek and Jech do not require there be a least element in the partial order which corresponds to trivial forcing, whereas for my context I find this a convenience, to fit with  trees always having a unique root. One can translate between the two different implementations by adding/removing the trivial forcing at the bottom.}
and follow with a more intuitive treatment.

\begin{definition}
Consider a tree $T$, possibly a proper class and suppose $\Pbb$ is a forcing partial order whose conditions are functions with domain a set-sized rooted subtree of $T$. For $S$ a rooted subtree of $T$, $s \in T$, and $p \in \Pbb$, make the following definitions:
\begin{itemize}
\item $\Pbb \rest S = \{ p \rest S : p \in \Pbb\}$ is the restriction of $\Pbb$ to $S$.
\item $\Pbb \rest \mathord{<}s = \Pbb \rest \{ t \in T : t < s \}$ is the restriction of $\Pbb$ to below $s$.
\item $\Pbb \rest \mathord{\le}s = \Pbb \rest \{ t \in T : t \le s \}$ is the restriction of $\Pbb$ up to $s$.
\item $p \rest s = p \rest (\Pbb \rest \mathord{<}s)$ is the restriction of $p$ to below $s$.
\end{itemize}
\end{definition}

\begin{definition}
Consider a tree $T$ and partial order $\Pbb$ as in the previous definition. Say that $\Pbb$ is a \emph{tree iteration} along $T$ if there is a $T$-indexed sequence $\seq{\dot \Qbb_s : s \in T}$ of names for \emph{iterands} satisfying the following.
\begin{itemize}
\item $\Qbb_\emptyset$ is trivial forcing.
\item For non-root $s \in T$, the name $\dot \Qbb_s$ is a $(\Pbb \rest \mathord{<}s)$-name.
\item For each $s \in T$, the forcing $\Pbb \rest \mathord{\le}s$ is the two-step iteration $(\Pbb \rest \mathord{<}s) * \dot \Qbb_s$.
\item If $p \in \Pbb$ and $s \in T$, then $p(s)$ is a name for a condition in $\dot \Qbb_s$. 
\item For $p,q \in \Pbb$, we have $p \le q$ if and only if $(i)$ $\dom p \supseteq \dom q$, $(ii)$ $p \rest s \le q \rest s$, and $(iii)$ $p \rest s \forces p(s) \le q(s)$, where we quantify over all appropriate $s$.
\item If $p \in \Pbb$ and $q \in \Pbb \rest S$, for some rooted subtree $S$, so that $(p \rest S) \le q$ then the condition $r$ which extends $q$ by setting $r(s) = p(s)$ for each $s \in \dom p \setminus S$ is in $\Pbb$.
\end{itemize}
A generic $\bar G$ for a tree iteration can be thought of as a sequence of generics $G_s$ for each iterand $\Qbb_s$.
\end{definition}

Observe that if $S$ is a rooted subtree of $T$ then $\Pbb \rest S$ is a complete suborder of $\Pbb$, and that a generic $\bar G$ for $\Pbb$ is interdefinable with the collection of its restrictions to rooted subtrees. If $S \subseteq T$ is a rooted subtree set $\bar G \rest S$ to be $\bar G \rest (\Pbb \rest S)$ and similarly for $\bar G \rest \mathord{<}s$ and $\bar G \rest \mathord{\le}s$.

\begin{observation}
Let $\Pbb$ be a tree iteration along $T$ and suppose $\bar G$ is a generic for $\Pbb$. If $s$ and $t$ are incomparable nodes in $T$ with infimum $r$, then $G_s$ and $G_t$ are mutually generic over $\Vrm[\bar G \rest \mathord{\le}r]$. More generally, if $S$ is a rooted subtree of $T$ and $s,t \in T \setminus S$ are incomparable nodes, then $G_s$ and $G_t$ are mutually generic over $\Vrm[\bar G \rest S]$. \qed
\end{observation}

This observation gives the basis for understanding how to inductively build up tree iterations: climbing up a branch corresponds to doing an ordinary iteration, and splitting at a node corresponds to product forcing at that stage. This is the approach I will take to define tree iterations.

Observe that tree iterations subsume both products and standard iterations; a product is a tree iteration along a tree consisting of only a root and its successors, while an iteration is a tree iteration along a non-branching tree. As a fact of particular use to us, a product of tree iterations is itself a tree iteration, with the underlying tree obtained by introducing a new root and making the old roots its immediate successors.

Groszek and Jech give calculations for chain condition and closure properties for generalized iterations. I state them below, specialized to the case of tree iterations.

\begin{lemma}[{\cite[Lemma 3]{groszek-jech1991}}] \label{lem:gj-closure}
Suppose $\Pbb$ is a tree iteration along a well-founded tree $T$ and that each iterand in $\Pbb$ is $\mathord{<}\kappa$-closed (respectively, $\mathord{<}\kappa$-distributive). Then $\Pbb$ is $\mathord{<}\kappa$-closed (respectively, $\mathord{<}\kappa$-distributive). 
\end{lemma}

\begin{lemma}[{\cite[Lemma 4]{groszek-jech1991}}] \label{lem:gj-cc}
Assume the $\GCH$ and let $\kappa < \lambda$ be regular cardinals. If $\Pbb$ is a tree iteration along a well-founded tree $T$ and whenever $S$ is a subtree of $T$ with cardinality $\mathord{<}\kappa$ we have $\card{\Pbb \rest S} < \lambda$, then $\Pbb$ has the $\lambda$-cc.
\end{lemma}

For this lemma, you can drop the assumption of the $\GCH$ by more carefully stating the cardinality assumptions. Groszek and Jech state it with the $\GCH$ assumption, and I left it in that form because the models in this paper will satisfy the $\GCH$. 
\smallskip

From these two lemmata we can derive what I will call the \emph{safety lemma} for coding. It states that you only add Cohen subsets to a cardinal you intend to, ensuring that different parts of the coding do not interfere with each other.

\begin{lemma}[Safety lemma] \label{lem:safety}
Assume the $\GCH$. Let $\Pbb$ is a tree iteration along a well-founded tree $T$ where each iterand in $\Pbb$ is a set-support iteration of $\Add^L$ forcings. For a cardinal $\alpha$ suppose that there is at most one iterand $\dot \Qbb_s$ in $\Pbb$ which adds a Cohen (over $\Lrm)$ subset to $\alpha$. Further suppose that every other iterand $\dot \Qbb_t$ is either $\mathord{<}\alpha$-distributive or has the $\alpha^+$-cc. 
Then, after forcing with $\Pbb$, any Cohen subset of $\alpha$ was added by the iterand $\dot \Qbb_s$. In particular, if no iterand adds a Cohen subset to $\alpha$, then $\Pbb$ adds no Cohen subset to $\alpha$.
\end{lemma}

\begin{proof}[Proof sketch]
If $\Pbb$ includes an iterand $\dot \Qbb_s$ which adds a Cohen subset to $\alpha$, then we can split $\Pbb$ into three pieces: first, the rooted subtree consisting of all stages not at or above $s$; second, the stage $s$ itself; third, the stages above $s$. If the first and third pieces, each of which has no iterand adding a Cohen subset to $\alpha$, do not add any Cohen subsets to $\alpha$, then we are done. So it suffices to prove the case where no iterand adds a Cohen subset to $\alpha$.

This is proved by an induction up $T$. A subtree consisting only of small---i.e. $\mathord{<}\alpha$-distributive---iterands will be $\mathord{<}\alpha$-distributive, by Lemma~\ref{lem:gj-closure}, and so not add any Cohens to $\alpha$. A subtree consisting only of big---i.e. with the $\alpha^+$-cc---iterands will have the $\alpha^+$-cc, and thus not add any Cohens to $\alpha$. We inductively break down $T$ into a tree of smaller subtrees, each of which consists either uniformly of small iterands or uniformly of big iterands. None of the subtree forcings add any Cohens to $\alpha$, and so inductively the entire forcing $\Pbb$ will add no Cohens to $\alpha$.
\end{proof}

With some general theory now at our disposal, let us discuss the specific coding forcings which will be used. To start, we want to force to code a set of ordinals $X$ into the $k$-th inner mantle. 
Ensuring we have enough space to code everything is straightforward. Using our favorite absolutely definable bijection ${}^{<\omega}\Ord \to \Ord$ we can partition a coding region $R$ into class-sized subregions $R_s$ for each $s \in {}^{< \omega}\Ord$, and if $R$ were absolutely definable then so would be the $R_s$. Our coding forcing will then be a tree iteration of Cohen coding forcings along the tree ${}^{\le k}\Ord$.

\begin{definition}
Consider a set $x$ of ordinals with $\lambda = \sup x$,\footnote{Recall that we assume without loss of generality that $\sup x \not \in x$.}
and a finite ordinal $k > 0$. Let $T = {}^{\le k}\Ord$. Then the \emph{$k$-height tree coding forcing} $\Tbb_k(x)$ is a tree iteration along $T$ defined as follows.
\begin{itemize}
\item The initial iterand is trivial forcing. Set $x_\emptyset = x$. Now proceed inductively upward along $T$ to define the rest of the forcing. 
\item After forcing stage $s \in T$ we have a set $x_s$ of ordinals to code for the next level of stages. We will do so $\Ord$ often, via a set-support product using Cohen coding. The iterand for stage $s \cat \seq\xi$ is Cohen coding of $x_s$. Call the code produced at stage $s \cat \seq\xi$ by $x_{s \cat \seq\xi}$, which we can think of as a set of ordinals.
\end{itemize}
Call the generic produced for the entire $\Tbb_k(x)$ forcing $\bar c$. We can think of this generic as a tree of sets of ordinals $x_s$, or we can think of $\bar c$ as a sequence of Cohens $c_\alpha \subseteq \alpha$ for an appropriate class of cardinals. By the Self-encoding Lemma \ref{lem:se} the entire sequence can be recovered from just having the $c_\alpha$'s by themselves, or from just having every  individual $x_s$.
\end{definition}

After forcing with $\Tbb_k(x)$, digging down the sequence of inner mantles corresponds to climbing down the tree of generics. To see why this is we will need to understand how $\Tbb_k(x)$ relates to its restriction to certain subtrees. 

If $T$ is a well-founded tree and $\ell \in \omega$, let $T_{-\ell}$ denote the subtree consisting of nodes whose distance to the closest leaf node is at least $\ell$. So 
\[
T = T_{-0} \supseteq T_{-1} \supseteq \cdots \supseteq T_{-\ell} \supseteq \cdots
\]
is a descending sequence of subtrees, where at each step we chop off all the leaf nodes. Given a tree iteration $\Pbb$ along $T$, this then gives a descending sequence of complete subforcings $\Pbb_{-\ell} = \Pbb \rest T_{-\ell}$ of $\Pbb$. To understand this sequence is it enough to analyze what happens at one step, as $\Pbb_{-(\ell+1)} = (\Pbb_{-\ell})_{-1}$.

Work in the intermediate extension by $\Pbb_{-1}$. For any leaf node $s \in T$, we have all the data needed to interpret the iterand $\dot \Qbb_s$. Note however, that the names $\dot \Qbb_s$ are not $\Pbb_{-1}$-names.\footnote{Except in the edge case where $T$ is linear at least up to the penultimate level.}
Rather, they are names in the subposet consisting of the path from the root to below $s$. So for our use case, $\Qbb_s$ is the Cohen coding forcing as computed in a certain inner model. To get from the intermediate extension to the full extension, we force with the (set support) product of these $\Qbb_s$'s. 
These iterands $\Qbb_s$'s are iterations of Cohen forcings, so we can think of their product as an iteration of Cohen forcings, except each that each stage is defined in $\Lrm$. This sort of iteration was previously studied by Reitz \cite{reitz2020}. Here is a definition specialized to our context.\footnote{Reitz's definition is more general, allowing iterands being of the form $\Add(\alpha,\lambda)$ (as defined in an inner model), rather than merely of the form $\Add(\alpha,1)$.}

\begin{definition}
A \emph{generalized Cohen iteration} is at iteration $\Pbb = \seq{\Pbb_\alpha,\dot \Qbb_\alpha : \alpha \in R}$ along a class $R$ of regular cardinals satisfying the following.
\begin{itemize}
\item $\dot \Qbb_\alpha$ is a name which is forced to be $\Add(\alpha,1)^{W_\alpha}$, where
\item $W_\alpha$ is (forced to be) an inner model of $\Vrm^{\Pbb_\alpha}$, one which contains the ground model $\Vrm$, and
\item $W_\alpha \models \card{\Add(\alpha,1)^V} = \card{\Add(\alpha,1)^{W_\alpha}}$.
\end{itemize}
\end{definition}

In the context of this article, we force over $\Lrm$ and $W_\alpha$ is always $\Lrm$, so the conditions are trivially satisfied. Note that this relies on it being, as it will be for us, that $\Lrm$ is correct about cardinals.

Reitz proved that generalized Cohen iterations have nice distributivity properties.

\begin{definition}
A class forcing $\Pbb$ is a \emph{progressively distributive iteration} if for every regular $\alpha$ you can factor $\Pbb$ as $\Pbb \cong \Pbb_\headsf * \dot \Pbb^\tailsf$ where $\Pbb_\headsf$ is a set and $\Pbb_\headsf \forces \dot \Pbb^\tailsf$ is $\mathord{<}\alpha$-distributive.
\end{definition}

\begin{theorem}[{Reitz \cite[Theorem 14]{reitz2020}}]
Any generalized Cohen iteration is a progressively distributive iteration. In particular, any generalized Cohen iteration preserves $\ZFC$.
\end{theorem}

When splitting an iteration consisting of Cohen coding forcings we may choose to keep all iterands in one Cohen coding forcing on the same side of the split, and thereby split this iteration into a \emph{product} $\Pbb_\headsf \times \Pbb^\tailsf$. Cohen coding forcing starting at $\alpha$ is $\mathord{<}\alpha$-distributive, and so Reitz's work yields that this forcing defined in an inner model which meets the requirements is $\mathord{<}\alpha$-distributive. So having the set-sized head-forcing $\Pbb_\headsf$ consisting of iterands from the Cohen coding for small enough starting points shows we have a \emph{progressively distributive product}.  That is, not only can we split the forcing into a set-sized head forcing and a sufficiently distributive tail forcing, but we can do these two forcings in either order.

We are (finally!) ready to compute some inner mantles. Let us begin with handling the situation with tree-like coding a single set of ordinals.

\begin{lemma} \label{lem:finite-tree-coding}
Force over $\Lrm[x]$ with $\Tbb = \Tbb_k(x)$ (using an appropriate coding region $R$) to get $\Lrm[x][\bar c]$. For $\ell \le k$ let $\bar c_{-\ell}$ be the restriction of $\bar c$ to $\Tbb_{-\ell}$. Then, for $\ell \le k$, we have $(\Mrm^\ell)^{\Lrm[x][\bar c]} = \Lrm[x][\bar c_{-\ell}]$. In particular, $(\Mrm^k)^{\Lrm[x][\bar c]} = \Lrm[x]$.
\end{lemma}

\begin{proof}
Work in $\Lrm[x][\bar c]$. We prove this by induction on $\ell$. The base case $\Mrm^0 = \Lrm[x][\bar c_{-0}] = \Lrm[x][\bar c]$ is trivial. 

$(\Mrm^{\ell+1} \supseteq \Lrm[x][\bar c_{-(\ell+1)}])$ Consider the piece of the generic $x_s$ coming from stage $s$ where $s \in T$ has$\length s < k - \ell$. By the inductive hypothesis, $x_{s \cat \seq \xi}$ is in $\Mrm^\ell = \Lrm[x][\bar c_{-\ell}]$ for every ordinal $\xi$. If $W \subseteq \Mrm^\ell$ is a ground then it contains $x_{s \cat \seq \xi}$ for some large enough $\xi$, whence $x_s \in W$. Since $W$ was arbitrary, we conclude $x_s \in \Mrm^{\Mrm^\ell}$. By the Self-encoding Lemma \ref{lem:se}, we may conclude the entire sequence $\bar c_{-(\ell+1)}$ is definable in $\Mrm^{\ell+1}$.

$(\Mrm^{\ell+1} \subseteq \Lrm[x][\bar c_{-(\ell+1)}])$ Consider a set $A \in \Lrm[x][\bar c_{-\ell}] \setminus \Lrm[x][\bar c_{-(\ell+1)}]$. Then $A$ was added by the $(k-\ell)$-th level of $\Tbb_k(X)$, call this forcing $\Sbb$. This forcing is a progressively distributive product. So we may split it as $\Sbb = \Sbb_\headsf \times \Sbb^\tailsf$ where $\Sbb_\headsf$ is a set and $\Sbb^\tailsf$ is $\card{A}^+$-distributive. If $G^\tailsf$ is the restriction of the generic $G \subseteq \Sbb$ to $\Sbb^\tailsf$, then $\Lrm[x][\bar c_{-(\ell+1)}][G^\tailsf]$ is a ground of $\Lrm[x][\bar c_{-\ell}]$ which does not contain $A$. Thus $A \not \in \Mrm^{\ell+1} = \Mrm^{\Lrm[x][\bar c_{-\ell}]}$.
\end{proof}

This lemma generalizes to products of tree-like codings.

\begin{lemma} \label{lem:prod-tree-coding}
Suppose $\Tbb$ is a product of tree-like codings $\Tbb_{k_i}(x_i)$ along disjoint coding regions, each of which is clean for coding. Think of $\Tbb$ as a tree iteration, where its nodes at level $n+1$ are the nodes of level $n$ along the $\Tbb_{k_i}(x_i,R_i)$'s. Let $\Lrm[\bar x][\bar c]$ be an extension by $\Tbb$, and set $\bar c_{-\ell}$ to be the restriction of $\bar c$ to $\Tbb_{-\ell}$. Then, for each $\ell < \omega$, we have that $(\Mrm^\ell)^{\Lrm[\bar x][\bar c]} = \Lrm[\bar x][\bar c_{-\ell}]$.
\end{lemma}

\begin{proof}[Proof sketch]
The same inductive argument of climbing down the tree as in the previous lemma applies. 
\end{proof}

\smallskip

To close off this section, I would like to briefly address the choice to use Cohen coding. For the results of this article this coding suffices. We are interested in constructing counterexample models, so it is no harm to use a coding method that only applies to a select class of models. And Cohen coding is particularly simple, which shortens some arguments. But the choice of coding is not essential, and one could use alternate methods. For instance, one could code sets by the pattern of where the $\GCH$ holds/fails, as was done in \cite{reitz-williams2019}.

\section{Failing to satisfy  \texorpdfstring{$\AC$}{AC}} \label{sec:choice}

In this section I construct a model whose $\omega$-th inner mantle is a definable class but fails to satisfy $\AC$, establishing Main Theorem~\ref{main:choice}. My strategy follows that of McAloon for the analogous result about $\HOD$ \cite{mcaloon1974}.

\begin{theorem}[An elaboration of Main Theorem~\ref{main:choice}]
There is a class forcing extension of $\Lrm$ in which $\Mrm^\omega$ is a definable class but $\Mrm^\omega \models$ ``there is no well-order of $\powerset(\omega)$''. 
\end{theorem}

\begin{proof}
Begin by forcing over $\Lrm$ to add a generic $A \subseteq \Add^\Lrm(\omega,\omega_1)$. To clarify, since the implementation details will matter later: Here $A$ is an $\omega \times \omega_1$ grid of $0$s and $1$s, where the rows of $A$ are the $\omega_1$ many Cohen reals added over $\Lrm$. For $k < \omega$ let $A^k = \{ (n,\alpha) \in A : n \ge k \}$ be the portion of $A$ from column $k$ rightward. 
Using our favorite absolutely definable pairing function we may think of each $A^k$ as a subset of $\omega_1$. In an abuse of notation, I will use $A^k$ to mean either the $\omega \times \omega_1$ binary grid and the subset of $\omega_1$ coding this grid, depending on context. We now want to force to code each $A^k$ into $\Mrm^k$ while ensuring that no $A^k$ is in $\Mrm^\omega$.

To do this, force with the full support product of the forcings $\Tbb_k(A^k)$, where these codings take place on disjoint regions. We have all of the cardinals $>\omega_1$ clean for coding, and we may split them into $\omega$ many subregions in an absolutely definable way. Call this product $\Tbb$, which we may think of as a tree iteration. Let $\bar c \subseteq \Tbb$ be generic over $\Lrm[A]$.

Work in $\Lrm[A][\bar c]$. By Lemma \ref{lem:prod-tree-coding} we can calculate the first $\omega$ many inner mantles. Namely, $\Mrm^k = \Lrm[\bar c_{-k}]$, where $\bar c_{-k}$ is the restriction of $\bar c$ to the nodes at most $k$ far from a leaf. In particular, $A^k \in \Mrm^k$. However, $A^k \not \in \Mrm^{k+1}$ because $A^k$ is not constructible from $A^{\ell}$ for $\ell > k$. 

It follows immediately from this computation of the inner mantles that the sequence of inner mantles $\seq{\Mrm^\eta : \eta \le \omega}$ is definable over $\Lrm[A][\bar c]$, and that in each inner mantle the tail sequence is definable.
By Lemma~\ref{lem:zf-is-free} we get that $\Mrm^\omega \models \ZF$. It remains only to see that $\Mrm^\omega$ does not have a well-order of its reals.
Observe that $\omega_1^{\Lrm[A][\bar c]} = \omega_1^\Lrm$, since adding $\bar c$ does not add any new reals. Thus, if $\Mrm^\omega$ has a well-order of its $\powerset(\omega)$ it must have one of ordertype $\omega_1$. In such a case, there would be $X \subseteq \omega_1$ so that $\powerset(\omega) \cap \Mrm^\omega \subseteq \Lrm[X]$. It therefore suffices to show there is no such $X$. 

Take any $X \subseteq \omega_1$ in $\Mrm^\omega$. Then $X \in \Lrm[A][\bar c]$, but we added $\bar c$ by coding high enough so as to not add new subsets of $\omega_1$. So already $X \in \Lrm[A]$. Because $\Add(\omega,\omega_1)$ has the ccc we get that $X$ was added by a countable piece of $A$. That is, there is an ordinal $\alpha < \omega_1$ so that $X \in \Lrm[A \rest \alpha]$ where $A \rest \alpha$ is the rows of $A$ below $\alpha$. Let $z$ be one of the rows of $A$ above $\alpha$. Then $z$ is Cohen-generic over $\Lrm[A \rest \alpha]$ and in particular $z \not \in \Lrm[X]$. But $z$ is in $\Mrm^\omega$ since each $\Mrm^k$ contains a tail of $z$ and $z$ can be recovered from a tail by prepending a finite binary sequence. We have seen that $z$ witnesses that $\powerset(\omega) \cap \Mrm^\omega \not \subseteq \Lrm[X]$, as desired to complete the proof of the theorem.
\end{proof}

Combining this theorem with known techniques we can see that this situation can happen in a model of $\Vrm = \HOD$.

\begin{corollary}
There is model of $\ZFC + \Vrm = \HOD$ in which $\Mrm^\omega$ is a definable class but fails to satisfy $\AC$.
\end{corollary}

\begin{proof}
Let $N$ be the extension of $\Lrm$ from the theorem. Class force further over $N$ to get a model $N[G] \models \Vrm = \HOD$ in which $\Mrm^{N[G]} = N$, as in \cite[Theorem 67]{fuchs-hamkins-reitz2015}. Then $(\Mrm^\omega)^{N[G]} = (\Mrm^\omega)^N$ is a definable class in $N[G]$ but fails to satisfy $\AC$.
\end{proof}

\section{Failing to be a definable class} \label{sec:definable}

In this section I  construct a model of $\ZFC$ whose $\omega$-th inner mantle fails to be a definable class, establishing Main Theorem~\ref{main:definable}. Following Harrington's analogous construction for the $\omega$-th $\HOD$,\footnote{See \cite[Section 7]{zadrozny1983}.}
the strategy will be as follows. We will start with a model of $\Vrm = \Lrm$. We will then force over this model to add certain Cohen subsets in such a way that which cardinals have Cohen subsets (generic over $\Lrm$) in $\Mrm^\omega$ codes the truth predicate for $\Lrm$. If our original model of $\Vrm = \Lrm$ is moreover a \emph{Paris model}, one whose ordinals are each definable without parameters,\footnote{Paris showed \cite{paris1973} that every consistent extension of $\ZF$ has a Paris model. Note that in particular the Shepherdson--Cohen \cite{shepherdson1953,cohen1963} minimum transitive model of $\ZF$ is a transitive Paris model of $\Vrm = \Lrm$.} 
then this will imply that $\Mrm^\omega$ cannot be a definable class. For if it were definable then we could define a surjection $\omega \to \Ord$, which would contradict the Replacement axiom schema.

To prove Main Theorem \ref{main:definable} we will need two additional tools beyond the coding forcings in Section~\ref{sec:coding}. 

The first, which I will call \emph{robust Cohen coding}, will not only ensure that the full sequence for the generic is recoverable from the individual Cohens, but will also ensure the same happens in every ground. To do this, rather than have $\omega$ stages in the iteration to code previous stages, have $\Ord$ many stages. In more detail, robust Cohen coding to code a set $x$, call it $\RCC(x)$, is a set-support $\Ord$-length iteration, which we will do along an absolutely definable coding region consisting of regular cardinals. Before stage $\xi+1$ in the iteration we have the sequence $\seq{c_i : i \le \xi}$ of generics up to stage $\xi$, which we may think of as a single set of ordinals $x_\xi$. For stage $\xi+1$ do simple Cohen coding forcing to code $x_\xi$ into the next unused block of the coding region. As with Cohen coding, we can compute in advance the blocks where each stage of the iteration will take place, even though we cannot compute in advance exactly which cardinals will be given Cohen generic subsets. It is straightforward to check that robust Cohen coding enjoys the self-encoding property of Cohen coding, and moreover that the entire generic can be recovered just from a tail of the individual Cohen subsets. In particular, any ground of a forcing by robust Cohen coding must contain the entire robust generic.

The second tool we need is a way to overwrite coding blocks. If $R$ is a class of cardinals used as a coding block, we can overwrite the coding pattern in $R$ by including a Cohen subset to every cardinal in $R$. Let $\Obb(R)$, the \emph{overwrite forcing on $R$} be the set-support product forcing which does this, namely
\[
\Obb(R) = \prod_{\alpha \in R} \Add^\Lrm(\alpha,1).
\]
It may be that the original codes living in the region $R$ are still definable---for example, if they were coded elsewhere---but in the absence of such this can be used to destroy coding information.

\begin{theorem} \label{thm:truth-in-momega}
There is a class-forcing extension of $\Lrm$ so that in the $\Mrm^\omega$ of this extension the truth predicate for $\Lrm$ is a definable class.
\end{theorem}

\begin{proof}
Fix an absolute one-to-one correspondence between formulae $\phi$ with parameters from $\Lrm$ and the class consisting of all cardinals $\kappa(\phi)$ with odd index, i.e. $\kappa(\phi) = \aleph_{{2\cdot\xi+1}}$ for some ordinal $\xi$. The class forcing we use will be a set-support product of forcings which add a Cohen subset $c_\phi$ to $\kappa(\phi)$ and follows this with a coding forcing which ensures $c_\phi$ is in $\Mrm^\omega$ if and only if $\Lrm \models \phi$. To do all this coding, we need to split up remaining cardinals into $\Ord$-sized coding regions for each $\phi$. Because $\Lrm$ has a $\Delta_1$ global well-order, this can be arranged. Each of these coding regions will in turn be split into yet smaller coding regions, one for each short enough sequence $\seq{x_1,\ldots,x_k}$ of sets from $\Lrm$. That too can be done in an absolute way. I leave the details of defining the coding regions in a way that ensures there is enough space for the following work to the interested reader. 

When defining the coding forcing for $c_\phi$, we must take care that the definition refers only to a bounded level of truth in $\Lrm$, to ensure that the entire product forcing is definable. Before jumping into the full details, let me illustrate the strategy by sketching what happens with a simple example. See Figure~\ref{fig:smoltriangle} for a pictorial representation of this example. 

\smallskip

Consider a $\Sigma_2$-formula $\phi$ in the form $\exists x \neg \exists y\  \psi(x,y)$ where $\psi$ is $\Sigma_0$ and possibly contains parameters from $\Lrm$. Consider the Cohen generic $c_\phi \subseteq \kappa(\phi)$. 
After adding this real, we force with a product, for each $x \in \Lrm$, of robust Cohen forcings to code $c_\phi$. Call $d_\phi(x)$ the robust generic corresponding to the choice of $x$, and $d_\phi(x;\xi)$ will denote the $\xi$-th set-sized piece of the generic, of which there are $\Ord$ many. 
Next, for each $d_\phi(x;\xi)$ force with the product of $\Tbb_1(d_\phi(x;\xi))$, producing a code $\bar d_\phi(x;\xi)$, and overwrite forcing for the region where $d_\phi(x;\xi)$ lives, producing an overwrite block $o_\phi(x;\xi)$ consisting of a Cohen subset for each cardinal in this region. Using a pairing function we may think of $o_\phi(x;\xi)$ as a single set of ordinals. Finally, for each $y \in \Lrm$ we query the $\Sigma_0$-truth predicate for $\Lrm$ and ask whether $\Lrm \models \psi(x,y)$. If so, follow the forcing to add $o_\phi(x;\xi)$ with $\Tbb_1(o_\phi(x;\xi))$, producing the code $\bar d_\phi(x,y;\xi)$. Else, if $\Lrm \not \models \psi(x,y)$, then do nothing.

Let us now analyze this and see that we have coded $c_\phi$ into $\Mrm^\omega$ if and only if $\Lrm \models \phi$. First, observe that $d_\phi(x;\xi)$ survives into $\Mrm$ for each $x \in \Lrm$ and each $\xi$. This is because any ground must contain a tail of $\bar d_\phi(x;\xi)$, and can thereby recover $d_\phi(x;\xi)$. Whether it survives to deeper inner mantles depends on whether $o_\phi(x;\xi)$ made it into $\Mrm$. Consider first the case where $\Lrm \models \exists x \neg \exists y\ \psi(x,y)$, as witnessed by a fixed $x \in \Lrm$. Then, there is no $y \in \Lrm$ so that $\Lrm \models \psi(x,y)$. Accordingly, we did not do any coding forcing to ensure $o_\phi(x;\xi)$ gets into $\Mrm$. So inside $\Mrm$ we can look at the region where $d_\phi(x;\xi+1)$ was coded and use the Cohen pattern there to define $d_\phi(x;\xi)$. So the entire robust generic $d_\phi(x)$ is a definable class in $\Mrm$, and hence $d_\phi(x) \in \Mrm^\omega$. This then implies that $c_\phi \in \Mrm^\omega$, as desired.

Consider next the other case that $\Lrm \not \models \exists x \neg \exists y\ \psi(x,y)$. Then, for each $x \in \Lrm$ there is $y \in \Lrm$ so that $\Lrm \models \psi(x,y)$. Fix $x \in \Lrm$ and and ordinal $\xi$. We coded $o_\phi(x;\xi)$ via some code $\bar d_\phi(x,y;\xi)$ and so $o_\phi(x;\xi)$ survives into $\Mrm$. Thus, in $\Mrm$ the coding region for each $d_\phi(x;\xi)$ is overwritten and so the robust generic $d_\phi(x)$ fails to be a definable class in $\Mrm$. No piece $d_\phi(x;\xi)$ will survive into $\Mrm^2$. Since this happens for every $x$, we get that $c_\phi \not \in \Mrm^\omega$.

Observe that the definition of which codes to include only used the $\Sigma_0$-truth predicate for $\Lrm$. We were able to code the $\Sigma_2$-truth predicate for $\Lrm$ into the $\omega$-th mantle of the extension by using a smaller fragment of the truth predicate for $\Lrm$. The task for the general case is then to see how we can do this for all $\Sigma_k$-formulae, reducing everything down to the $\Sigma_0$-truth predicate for $\Lrm$. We will follow the idea of $\Sigma_2$-case, but with more layers of coding.

\begin{figure}
\begin{center}

\newcommand\cod{D}
\begin{tikzpicture}[n/.style={align=center,scale=.7},>=stealth]

\node[n] (x)   at   (0,0) {{\Large $\Add(\kappa_\phi,1)$}};
\node[n] (d1)  at   (0,2) {{\Large $\mathbb{RCC}$} \\ for each $x$};
\node[n] (o1)  at   (3,2) {{\Large $\Obb$} \\ for each $x;\xi$};
\node[n] (dk)  at (3,4) {{\Large $\Tbb_1$} \\ for each $x,y;\xi$ \\ so that $\Lrm \models \psi(x,y)$};

\node[n] (t1)  at  (0,4) {{\Large $\Tbb_{1}$} \\ for each $x;\xi$};

\draw[->] (t1)  -- (d1);
\draw[->] (d1)  -- (x);
\draw[->] (dk)  -- (o1);

\draw[->,dashdotted] (o1)  -- (d1);

\end{tikzpicture}
\end{center}
\caption{Triangle coding for the simple $\Sigma_2$ example. Each node in the diagram says what forcing is done at each stage and the parameters determining how many copies of the forcing are done at each stage.}
\label{fig:smoltriangle}
\end{figure}
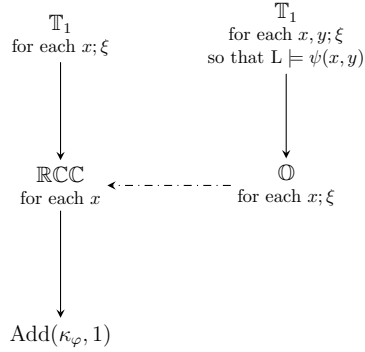

\begin{definition*}[Triangle coding]
Let me now describe the general case, a tree iteration I will call \emph{triangle coding}. Fix a $\Sigma_k$-formula $\phi$ with parameters from $\Lrm$, where $k > 0$ is even.\footnote{A small modification would handle the odd case, but it suffices to consider only the even case and it slightly simplifies the presentation.}
Write $\phi$ in the form
\[
\exists x_1 \neg \exists x_2 \exists \cdots \neg \exists x_k\ \psi(x_1, \ldots, x_k),
\]
where $\psi$ is $\Sigma_0$.

to inductively describe the triangle coding forcing. I will organize it into innings, which consist of stages on multiple levels of the tree, starting with Inning $-2$ and going to Inning $k$. See Figure~\ref{fig:triangle} for a pictorial representation of the coding pattern.
\begin{itemize}
\item (Inning $-2$) The root stage is trivial forcing.
\item (Inning $-1$) The first nontrivial level of the of the forcing consists of a single stage, adding the Cohen generic $c_\phi \subseteq \kappa(\phi)$. 
\item (Inning $0$) The second level of the forcing, consists of one stage for each $x_1 \in \Lrm$. For each stage do robust Cohen coding to code $c_\phi$, producing a code $d_\phi(x_1)$. This code is an $\Ord$-length sequence of generics, each coding the previous generics. Call the set-sized pieces $d_\phi(x_1;\xi_1)$, for $\xi_1 \in \Ord$.
\item (Inning $1$) We continue building upward with two subtrees above this stage, which I will call the above subtree and the right subtree. The above subtree consists of the set-support product of $\Tbb_{k-1}(d_\phi(x_1;\xi_1))$, for $\xi_1 \in \Ord$, to ensure that the codes $d_\phi(x_1;\xi_1)$ survive into $\Mrm^{k-1}$. Whether they survive deeper depends on what happens in the right subtree. Call the generics from these $\Tbb_{k-1}(d_\phi(x_1;\xi_1))$ by $\bar d_\phi(x_1;\xi_1)$.

The immediate next level in the right subtree consists of overwrite forcings in the regions where the $d_\phi(x_1;\xi_1)$ are coded. As usual, we take a set-support product to extend the tree iteration upward to this next level. Call the resulting generics by $o_\phi(x_1;\xi_1)$, and we will think of each as a single set of ordinals. We inductively continue building up the right subtree in the subsequent innings.
\item (Inning $\ell$, for $1 < \ell < k$) Immediately prior to this inning we did overwrite forcings to produce generics $o_\phi(x_1,\ldots, x_{\ell-1};\xi_1,\ldots,\xi_{\ell-1})$. We need to say how we continue the tree iteration past these stages. For the immediate next stages, for each $x_\ell \in \Lrm$, do $\Tbb_1(o_\phi(x_1,\ldots, x_{\ell-1};\xi_1,\ldots,\xi_{\ell-1}))$, to produce a code $d_\phi(x_1,\ldots,x_\ell;\xi_1,\ldots,\xi_{\ell-1})$. Recall that this forcing is an $\Ord$-sized product of forcings to code the overwrite generic into the Cohen pattern. We can break this product up into set-sized pieces, splitting each class-sized code into set-sized pieces $d_\phi(x_1,\ldots,x_\ell;\xi_1,\ldots,\xi_\ell)$ indexed by $\xi_\ell \in \Ord$.

We continue the tree iteration beyond each $\Tbb_1$, with two subtrees beyond each stage. The above subtree for the stage indexed by $\seq{x_1,\ldots,x_\ell;\xi_1,\ldots,\xi_{\ell-1}}$ consists of a set-support product of $\Tbb_{k-\ell}(d_\phi(x_1,\ldots,x_\ell;\xi_1,\ldots,\xi_\ell))$ for each $\xi_\ell \in \Ord$. This part of the tree iteration will ensure that the codes $d_\phi(x_1,\ldots,x_\ell;\xi_1,\ldots,\xi_\ell)$ survive into $\Mrm^{k-\ell}$. Whether they survive deeper depends on what happens in the right subtree. Call the generics from these tree-like codings by $\bar d_\phi(x_1,\ldots,x_\ell;\xi_1,\ldots,\xi_\ell)$.

The immediate next level in the right subtree consists of overwrite forcings in the regions where the $d_\phi(x_1,\ldots,x_\ell;\xi_1,\ldots,\xi_\ell)$ are coded. As usual, we take a set-support product to extend the tree iteration upward to this next level. Call the resulting generics by $o_\phi(x_1,\ldots,x_\ell;\xi_1,\ldots,\xi_\ell)$, and we will think of each as a single set of ordinals. We inductively continue building up the right subtree in the subsequent innings.
\item (Inning $k$) We are almost done. In right subtree piece of Inning $k-1$ we produced overwrite generics $o_\phi(x_1,\ldots,x_{k-1};\xi_1,\ldots,\xi_{k-1})$, and now we have to decide whether to force to code them into $\Mrm^1$. To make this decision we will query the $\Sigma_0$-truth predicate for $\Lrm$, once for each $x_k \in \Lrm$. If $\Lrm \models \psi(x_1,\ldots, x_k)$ then we will have a stage consisting of $\Tbb_1(o_\phi(x_1,\ldots,x_{k-1};\xi_1,\ldots,\xi_{k-1}))$, producing a code $d_\phi(x_1,\ldots,x_k;\xi_1,\ldots,\xi_{k-1})$. Else, if $\Lrm \not \models \psi(x_1,\ldots, x_k)$, then we include no stage corresponding to $x_k$. As usual, these stages extend the previous one via a set-support product.
\end{itemize}
After these $k$ innings we have a tree iteration along a well-founded tree. As we will see, this triangle coding forcing ensures that $c_\phi$ survives into $\Mrm^\omega$ if and only if $\Lrm \models \phi$.
\end{definition*}
\smallskip

To get the full class forcing we will use, take the set-support product of the triangle codings for each $\phi$. This product of tree iterations, call it $\Pbb$, is itself a tree iteration along a well-founded tree.

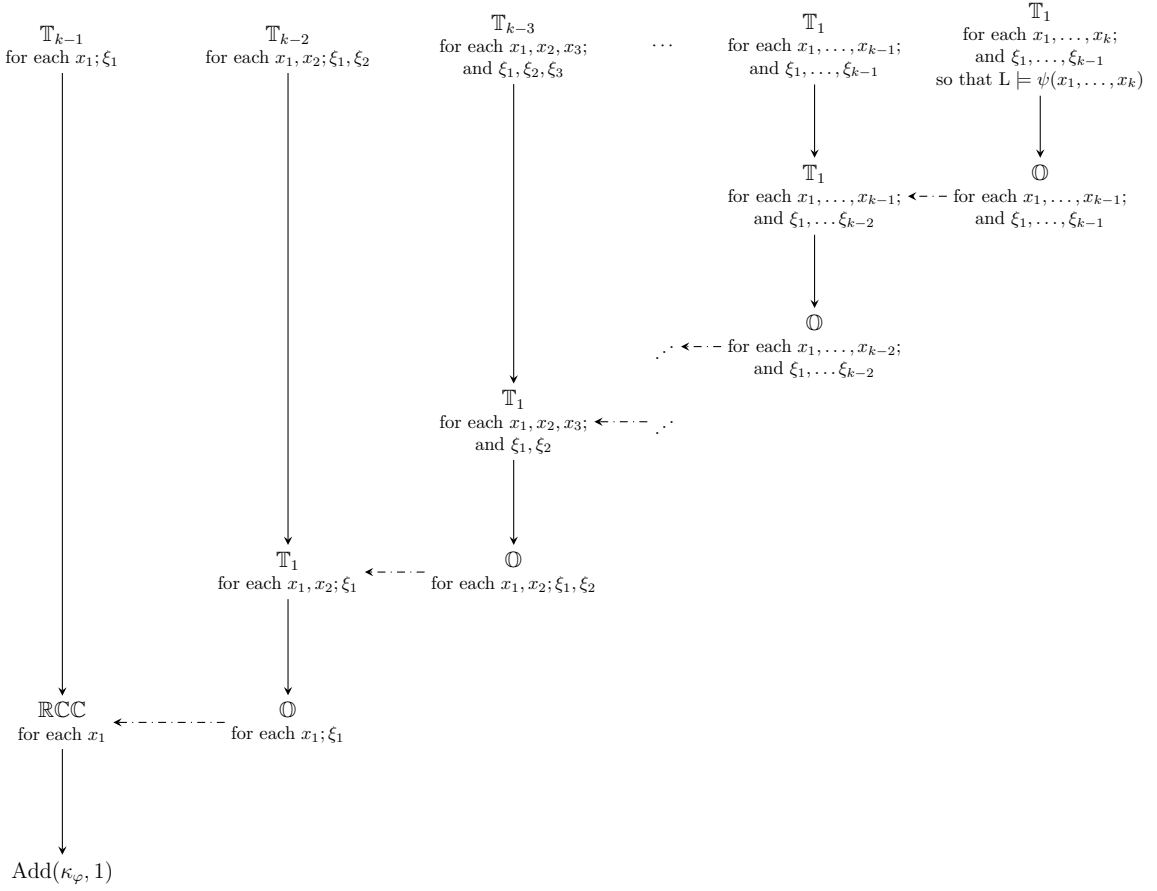
\begin{figure}
\begin{center}

\newcommand\cod{D}
\begin{tikzpicture}[n/.style={align=center,scale=.7},>=stealth]

\node[n] (x)   at   (0,0) {{\Large $\Add(\kappa_\phi,1)$}};
\node[n] (d1)  at   (0,2) {{\Large $\mathbb{RCC}$} \\ for each $x_1$};
\node[n] (o1)  at   (3,2) {{\Large $\Obb$} \\ for each $x_1;\xi_1$};
\node[n] (d2)  at   (3,4) {{\Large $\Tbb_1$} \\ for each $x_1,x_2;\xi_1$};
\node[n] (o2)  at   (6,4) {{\Large $\Obb$} \\ for each $x_1,x_2;\xi_1,\xi_2$};
\node[n] (d3)  at   (6,6) {{\Large $\Tbb_1$} \\ for each $x_1,x_2,x_3;$ \\ and $\xi_1,\xi_2$};
\node[n] (o-2) at  (10,7) {{\Large $\Obb$} \\ for each $x_1,\ldots,x_{k-2};$ \\ and $\xi_1,\ldots\xi_{k-2}$};
\node[n] (d-1) at  (10,9) {{\Large $\Tbb_1$} \\ for each $x_1,\ldots,x_{k-1};$ \\ and $\xi_1,\ldots\xi_{k-2}$};
\node[n] (o-1) at  (13,9) {{\Large $\Obb$} \\ for each $x_1,\ldots, x_{k-1};$ \\ and $\xi_1,\ldots,\xi_{k-1}$};
\node[n] (dk)  at (13,11) {{\Large $\Tbb_1$} \\ for each $x_1,\ldots, x_k;$ \\ and $\xi_1,\ldots,\xi_{k-1}$ \\ so that $\Lrm \models \psi(x_1,\ldots,x_k)$};

\node[n] (t1)  at  (0,11) {{\Large $\Tbb_{k-1}$} \\ for each $x_1;\xi_1$};
\node[n] (t2)  at  (3,11) {{\Large $\Tbb_{k-2}$} \\ for each $x_1,x_2;\xi_1,\xi_2$};
\node[n] (t3)  at  (6,11) {{\Large $\Tbb_{k-3}$} \\ for each $x_1,x_2,x_3;$ \\ and $\xi_1,\xi_2,\xi_3$};
\node[n] (t-1) at (10,11) {{\Large $\Tbb_1$} \\ for each $x_1,\ldots,x_{k-1};$ \\ and $\xi_1,\ldots,\xi_{k-1}$};

\node[n] (eb1) at  (8,6) {$\iddots$};
\node[n] (eb2) at  (8,7) {$\iddots$};
\node[n] (et)  at (8,11) {$\cdots$};

\draw[->] (t1)  -- (d1);
\draw[->] (d1)  -- (x);
\draw[->] (t2)  -- (d2);
\draw[->] (d2)  -- (o1);
\draw[->] (t3)  -- (d3);
\draw[->] (d3)  -- (o2);
\draw[->] (t-1) -- (d-1);
\draw[->] (d-1) -- (o-2);
\draw[->] (dk)  -- (o-1);

\draw[->,dashdotted] (o1)  -- (d1);
\draw[->,dashdotted] (o2)  -- (d2);
\draw[->,dashdotted] (eb1) -- (d3);
\draw[->,dashdotted] (o-2) -- (eb2);
\draw[->,dashdotted] (o-1) -- (d-1);

\end{tikzpicture}
\end{center}
\caption{Triangle coding to code $c_\phi$ into $\Mrm^\omega$ if and only if $\Lrm \models \phi$. Each node in the diagram says what forcing is done at each stage and the parameters determining how many copies of the forcing are done at each stage. Solid arrows denote what is being coded, and dashed arrows denote what is being overwritten.}
\label{fig:triangle}
\end{figure}

Work in $\Lrm[G]$, the extension by forcing with $\Pbb$. Let us analyze the sequence of inner mantles, which I will organize as a series of lemmata. The Safety Lemma~\ref{lem:safety} implies that to understand where $c_\phi$ is coded it suffices to look only at the iterand of the product $\Pbb$ corresponding to $\phi$.

\begin{sublemma} \label{slem:dphi}
For each $\Sigma_k$-formula $\phi$, sequence $\seq{x_1,\ldots,x_\ell}$ from $\Lrm$, and sequence $\seq{\xi_1,\ldots,\xi_\ell}$ of ordinals, where $\ell \le k$, the code $d_\phi(x_1,\ldots,x_\ell;\xi_1,\ldots,\xi_\ell) \in \Mrm^{k-\ell}$.
\end{sublemma}

\begin{proof}
This is because we did $(k-\ell)$-height tree-like coding of the code $d_\phi(\cdots)$ in a region that was never overwritten. Like in Lemma~\ref{lem:finite-tree-coding}, the portion of the code which survives into $\Mrm^i$ is precisely that corresponding to nodes distance $\ge i$ from the top. After $k-\ell$ many steps all that remains is the set $d_\phi(x_1,\ldots,x_\ell;\xi_1,\ldots,\xi_\ell)$ being coded by this tree-like coding.
\end{proof}

Next, let us see the role of the overwrite generics. We start with the penultimate inning, and work our way down the tree.

\begin{sublemma} \label{slem:sigma0}
For each $\Sigma_k$-formula $\phi$, each sequence $\seq{x_1,\ldots,x_{k-1}}$ from $\Lrm$ and each sequence $\seq{\xi_1,\ldots,\xi_{k-1}}$ of ordinals, the overwrite generic $o_\phi(x_1,\ldots,x_{k-1};\xi_1,\ldots,\xi_{k-1})$ is in $\Mrm^1$ if and only there is $x_k \in \Lrm$ so that $\Lrm \models \psi(x_1,\ldots,x_k)$. 
\end{sublemma}

\begin{proof}
$(\Rightarrow)$ Fix $x_k$ so that $\Lrm \models \psi(x_1,\ldots,x_k)$. Then, the code $d_\phi(x_1,\ldots,x_k;\xi_1,\ldots,\xi_{k-1})$ ensures that the overwrite generic gets into $\Mrm^1$. Specifically, any ground will contain a tail of the code $d_\phi(x_1,\ldots,x_k;\xi_1,\ldots,\xi_{k-1})$, and thus any ground can construct the overwrite generic by looking at the Cohen pattern in the region where $d_\phi(x_1,\ldots,x_k;\xi_1,\ldots,\xi_{k-1})$ is coded.

$(\Leftarrow)$ Suppose $\Lrm \not \models \psi(x_1,\ldots,x_k)$ for every $x_k \in \Lrm$. Then in the $\Pbb$ forcing there are no stages beyond the stage which added the overwrite generic. Let $G^-$ be the portion of the generic which excludes this overwrite forcing stage. Then $V[G^-]$ is a ground which does not contain the overwrite generic. (Namely, to get to $V[G]$ from $V[G^-]$ you force with the overwrite forcing, which is a set-sized forcing.) So the overwrite generic is not in the mantle.
\end{proof}

\begin{sublemma} \label{slem:succ1}
For each $\Sigma_k$-formula $\phi$, sequence $\seq{x_1,\ldots,x_\ell}$ from $\Lrm$, and sequence $\seq{\xi_1,\ldots,\xi_\ell}$ of ordinals, where $\ell < k-1$, the overwrite generic $o_\phi(x_1,\ldots,x_\ell;\xi_1,\ldots,\xi_\ell)$ is in $\Mrm^{k-\ell+1}$ if and only every overwrite generic $o_\phi(x_1,\ldots,x_{\ell+1};\xi_1,\ldots,\xi_{\ell+1})$, across all $x_{\ell+1}$, is \emph{not} in $\Mrm^{k-\ell}$. 
\end{sublemma}

\begin{proof}
From Lemma~\ref{slem:dphi} we know that the codes $d_\phi(x_1,\ldots,x_\ell;\xi_1,\ldots,\xi_\ell)$ survive into $\Mrm^{k-\ell}$.
Consider the case where overwrite generic from the right $o_\phi(x_1,\ldots,x_{\ell+1};\xi_1,\ldots,\xi_{\ell+1})$ is not in $\Mrm^{k-\ell}$. Then, any ground of $\Mrm^{k-\ell}$ can see a tail of the Cohen coding pattern in the region where the code $d_\phi(x_1,\ldots,x_\ell;\xi_1,\ldots,\xi_\ell)$ lives, and thereby recover the overwrite generic which is coded by $d_\phi(x_1,\ldots,x_\ell;\xi_1,\ldots,\xi_\ell)$, namely $o_\phi(x_1,\ldots,x_\ell;\xi_1,\ldots,\xi_\ell)$. Thus, this overwrite generic must get into $\Mrm^{k-\ell+1}$. 

Consider next the case where the overwrite generic from the right $o_\phi(x_1,\ldots,x_{\ell+1};\xi_1,\ldots,\xi_{\ell+1})$ is in $\Mrm^{k-\ell}$. I claim that in this case, $\Mrm^{k-\ell}$ has a ground which is the forcing to add the below overwrite generic $o_\phi(x_1,\ldots,x_\ell;\xi_1,\ldots,\xi_\ell)$. Each cardinal $\alpha$ in the region where $d_\phi(x_1,\ldots,x_\ell;\xi_1,\ldots,\xi_\ell)$ lives either had one or two Cohen (over $\Lrm$) generics added. Since $\Add^\Lrm(\alpha,1)$ is forcing equivalent to $\Add^\Lrm(\alpha,2)$, this amounts to having a single Cohen generic for each $\alpha$. The missing information to recover $d_\phi(x_1,\ldots,x_\ell;\xi_1,\ldots,\xi_\ell)$ is its domain, which is interdefinable with the below overwrite generic $o_\phi(x_1,\ldots,x_\ell;\xi_1,\ldots,\xi_\ell)$.
\end{proof}

\begin{sublemma} \label{slem:succ2}
For each $\Sigma_k$-formula $\phi$, each $x_1 \in \Lrm$, the robust Cohen code $d_\phi(x_1)$ survives into $\Mrm^k$ if and only if the overwrite generic $o_\phi(x_1;\xi_1)$ is \emph{not} in $\Mrm^{k-1}$ for every ordinal $\xi_1$. 
\end{sublemma}

\begin{proof}
This is the similar the previous lemma, except I stated it separately because $d_\phi(x_1;\xi_1)$ codes something different than the $d_\phi(\cdots)$ codes from later innings. The same argument goes through, except that we use the properties of robust Cohen coding to know that the entire generic $d_\phi(x_1)$ is definable in $\Mrm^{k}$ from knowing its definable in $\Mrm^{k-1}$ by looking at the non-overwritten region where it lives. Note that we use that the choice of $\xi_1$ does not matter; inductively applying Lemma~\ref{slem:succ1} yields that if one choice of $\xi_1$ works then they all do.
\end{proof}

\begin{sublemma}
For each $\Sigma_k$-formula $\phi$, the Cohen generic $c_\phi \subseteq \kappa_\phi$ survives into $\Mrm^\omega$ if and only if $\Lrm \models \phi$.
\end{sublemma}

\begin{proof}
This follows from the previous two lemmata by an induction on the quantifiers in $\phi$, working from the innermost out. For each $0 < \ell \le k$ let $\phi_\ell$ be obtained from $\psi_0$ by prepending all but the outermost $\ell$ quantifier blocks from $\phi$. So $\phi$ can be written as $\exists x_1\ \phi_1$, or $\exists x_1 \neg \exists x_2\ \phi_2$, or so on.

The base case, that $o_\phi(x_1,\ldots,x_{k-1};\xi_1,\ldots,\xi_{k-1})$ is in $\Mrm^1$ if only if $\Lrm \models \exists x_k\ \psi(x_1,\ldots,x_k)$, is Lemma~\ref{slem:sigma0}. Break the inductive step, from $\ell$ to $\ell-1$ up into even and odd cases. 
\begin{itemize}
\item ($\ell$ is even) Assume that the overwrite generic $o_\phi(x_1,\ldots,x_{\ell-1};\xi_1,\ldots,\xi_{\ell-1})$ is in $\Mrm^{k-\ell}$ if and only if $\Lrm \models \exists x_\ell\ \phi_\ell(x_1,\ldots,x_\ell)$. Then, by either Lemma~\ref{slem:succ1}, a below overwrite generic $o_\phi(x_1,\ldots,x_{\ell-2};\xi_1,\ldots,\xi_{\ell-2})$ survives into $\Mrm{k-\ell+1}$ if and only if there is \emph{no} $x_\ell$ so that $\Lrm \models \phi_\ell(x_1,\ldots, x_\ell)$. So some below overwrite generic survives if and only if $\Lrm \models \exists x_{\ell-1}\neg \exists x_\ell\ \phi_\ell(x_1,\ldots,x_\ell)$. But this last formula is just $\exists x_{\ell-1}\ \phi_{\ell-1}(x_1,\ldots,x_{\ell-1})$.

\item ($\ell$ is odd) This is handled similarly to the even case, except swap which quantifier is negated. Note that for the $\ell = 1$ case we use Lemma~\ref{slem:succ2} instead of Lemma~\ref{slem:succ1}.
\end{itemize}
So after step $\ell=1$ we have that, for each $x_1 \in \Lrm$ that the robust Cohen code $d_\phi(x_1)$ survives into $\Mrm^k$ if and only if $\Lrm \models \neg \exists x_2\ \phi_2(x_1,x_2)$. But then $c_\phi \in \Mrm^\omega$ if and only if $\Lrm \models \exists x_1\neg \exists x_2\ \phi_2(x_1,x_2)$, i.e. $\Lrm \models \phi$, as desired.
\end{proof}

As a consequence of this final lemma, $\Mrm^\omega$ can define the truth predicate for $L$. Namely, define it as the class of all $\kappa_\phi$ which contain a Cohen generic subset. The lemma tells us $c_\phi \in \Mrm$ if and only $\phi$ is true in $\Lrm$ and the Safety Lemma~\ref{lem:safety} tells us that in case $\phi$ is not true in $\Lrm$ then there are no Cohen generic subsets of $\kappa_\phi$.
\end{proof}

\begin{corollary}[An elaboration of Main Theorem~\ref{main:definable}]
Any Paris model of $\ZFC + \Vrm = \Lrm$ can be extended to a model of $\ZFC$ whose $\Mrm^\omega$ is not a definable class.
\end{corollary}

\begin{proof}
Let $L$ be a Paris model of $\ZFC + \Vrm = \Lrm$. Let $L[G]$ be the class forcing extension as in Theorem~\ref{thm:truth-in-momega} so that the truth predicate for $L$ is a definable class in $(\Mrm^\omega)^{L[G]}$. Suppose toward a contradiction that $\Mrm^\omega$ is a definable class in $L[G]$. Then, the truth predicate for $L$ is a definable class in $L[G]$. Because every ordinal in $L$ is definable without parameters from this truth predicate we can extract a surjection $\omega \to \Ord^L$, by mapping (G\"odel numbers of) formulae defining an ordinal to the ordinal they define in $L$. But then in $L[G]$ we have a definable surjection $\omega \to \Ord^L$, so by Replacement $\Ord^L$ is a set in $L[G]$, which is impossible.
\end{proof}

Let me also mention what this construction gives for $\omega$-nonstandard models.

\begin{corollary}
Suppose $L$ is an $\omega$-nonstandard Paris model of $\ZFC + \Vrm = \Lrm$. Let $L[G]$ be the class forcing extension from Theorem~\ref{thm:truth-in-momega}. Then in $L[G]$ the sequence of inner mantles does not stabilize at any standard index $k$ but no inner mantle $\Mrm^e$ of nonstandard index $e$ is definable.
\end{corollary}

\begin{proof}
Work in $L[G]$. That the sequence of inner mantles does not stabilize at a finite index follows from the construction---we always code some information into $\Mrm^{k}$ which was not in $\Mrm^{k+1}$. Suppose toward a contradiction that $\Mrm^e$ were definable. Inside $\Mrm^e$ consider the class $T$ of all $\phi \in \Lrm$ so that $c_\phi$ has a Cohen generic (over $\Lrm$) subset. Externally to the model, we see from the work in the theorem that $T$ restricted to standard $\phi$ is the truth predicate for $L$. ($T$ will also contain some $\phi$ with nonstandard G\"odel code.) So in $\Mrm^e$ we can define a cofinal map $\omega \to \Ord$ by sending $n$ to an ordinal $\xi$ if $n$ is the G\"odel code for a formula $\phi(x)$ so the formulae ``$\phi(\xi)$'' and ``$\exists!x\ \phi(x)$'' are  in $T$, and sending $n$ to $0$ otherwise. But then $\Mrm^e$ fails to satisfy Replacement, contradicting that successor stages in the sequence of inner mantles always satisfy full $\ZFC$.
\end{proof}

We can make the cut of indices of definable versus undefinable inner mantles occur at any $\Zbb$-block in $\omega^L$, not just the standard cut.

\begin{corollary}
Suppose $L$ is an $\omega$-nonstandard Paris model of $\ZFC + \Vrm = \Lrm$. Fix a nonstandard $e \in \omega^L$. Then $L$ has a class forcing extension $L[G][H]$ in which the indices $f$ for definable inner mantles are precisely those $f \in \omega^L$ which are $\le e + n$ for some standard $n$.
\end{corollary}

\begin{proof}
Get $L[G]$ using the forcing from the theorem, and then further force with $\Mbb(e)$ from \cite{reitz-williams2019}, the class forcing which makes the ground model the $e$-th inner mantle of the extension (with the sequence being strictly decreasing up to $e$.
\end{proof}

\section{Future work}

It is natural to ask whether there is anything special about $\omega$, and I conjecture that there is not.

\begin{conjecture}
For any limit ordinal $\gamma$, it is consistent that the sequence of inner mantles does not stabilize before $\gamma$ and that $\Mrm^\gamma$ is a definable inner model of $\ZF + \neg \AC$.
\end{conjecture}

\begin{conjecture}
For any limit ordinal $\gamma$, it is consistent that the sequence of inner mantles does not stabilize before $\gamma$ and that $\Mrm^\gamma$ fails to be a definable class.
\end{conjecture}

\bibliographystyle{alpha}
\bibliography{ref}

\begin{thebibliography}{FHR15}

\bibitem[Coh63]{cohen1963}
Paul Cohen.
\newblock A minimal model for set theory.
\newblock {\em Bulletin of the American Mathematical Society}, 69(4):537--540,
  1963.

\bibitem[FHR15]{fuchs-hamkins-reitz2015}
Gunter Fuchs, Joel~David Hamkins, and Jonas Reitz.
\newblock Set-theoretic geology.
\newblock {\em Annals of Pure and Applied Logic}, 166(4):464--501, 2015.

\bibitem[GJ91]{groszek-jech1991}
M.~Groszek and T.~Jech.
\newblock Generalized iteration of forcing.
\newblock {\em Transactions of the American Mathematical Society}, 324:1--26,
  1991.

\bibitem[Lav07]{laver2007}
Richard Laver.
\newblock Certain very large cardinals are not created in small forcing
  extensions.
\newblock {\em Annals of Pure and Applied Logic}, 149(1):1--6, 2007.

\bibitem[McA74]{mcaloon1974}
K.~McAloon.
\newblock On the sequence of models $\mathrm{HOD}_n$.
\newblock {\em Fundamenta Mathematicae}, 82:85--93, 1974.

\bibitem[Par73]{paris1973}
J.~Paris.
\newblock {\em Minimal models of {ZF}}, pages 327--331.
\newblock Leeds University Press, 1973.

\bibitem[Rei20]{reitz2020}
Jonas Reitz.
\newblock Cohen forcing and inner models.
\newblock {\em Mathematical Logic Quarterly}, 66(1):65--72, 2020.

\bibitem[RW19]{reitz-williams2019}
Jonas Reitz and Kameryn~J. Williams.
\newblock Inner mantles and iterated hod.
\newblock {\em Mathematical Logic Quarterly}, 65(4):498--510, 2019.

\bibitem[She53]{shepherdson1953}
J.~C. Shepherdson.
\newblock Inner models for set theory--part {III}.
\newblock {\em The Journal of Symbolic Logic}, 18(2):145--167, 1953.

\bibitem[Usu17]{usuba2017}
Toshimichi Usuba.
\newblock The downward directed grounds hypothesis and very large cardinals.
\newblock {\em Journal of Mathematical Logic}, 17(02):1750009, 2017.

\bibitem[Woo11]{woodin2011b}
W.~Hugh Woodin.
\newblock The continuum hypothesis, the generic-multiverse of sets, and the
  {$\Omega$} conjecture.
\newblock In Juliette Kennedy and Roman Kossak, editors, {\em Set Theory,
  Arithmetic, and Foundations of Mathematics}, pages 13--42. Cambridge
  University Press, 2011.
\newblock Cambridge Books Online.

\bibitem[Zad83]{zadrozny1983}
W{\l}odzimierz Zadro{\.z}ny.
\newblock Iterating ordinal definability.
\newblock {\em Annals of Pure and Applied Logic}, 24(3):263--310, 1983.

\end{thebibliography}

\end{document}